\documentclass[11pt,fleqn]{article}
\RequirePackage[T1]{fontenc}
\RequirePackage{amsthm,amsmath}
\RequirePackage[round]{natbib}
\usepackage[english]{babel}
\RequirePackage[colorlinks,citecolor=blue,urlcolor=blue]{hyperref}
\usepackage{enumerate}
\usepackage{caption}
\usepackage{amsmath,amsfonts,amssymb,euscript,mathrsfs}
\usepackage{graphicx}
\usepackage{fullpage}
\usepackage{amsmath}
\usepackage{dsfont}
\makeatletter
\def\captionof#1#2{{\def\@captype{#1}#2}}
\makeatother
\def\1{\mbox{\bf 1}}
\def\R{\mathbb{R}}

\def\N{\mathbb{N}}
\def\P{\mathbb{P}}
\def\E{\mathbb{E}}

\def\R{\mathbb{R}}

\def\Z{\mathbb{Z}}
\def\G{\mathbb{G}}

\def\c{\mbox{Cov}}

\newtheorem{theo}{Theorem}
\newtheorem{lem}{Lemma}
\newtheorem{prop}{Proposition}

\newtheorem{cor}{Corollary}
\newtheorem{Def/Prop}{Definition-Proposition}

\newcounter{exos}
\renewcommand\theexos{\arabic{exos}}

\newcounter{prob}
\renewcommand\theprob{\arabic{prob}}

\begin{document}
\author{Lionel Truquet \footnote{CREST-ENSAI, UMR CNRS 9194, Campus de Ker-Lann, rue Blaise Pascal, BP 37203, 35172 Bruz cedex, France. {\it Email: lionel.truquet@ensai.fr}}}

\title{Strong mixing properties of discrete-valued time series with exogenous covariates}
\date{}
\maketitle

\begin{abstract}
We derive strong mixing conditions for many existing discrete-valued time series models that include exogenous covariates in the dynamic.
Our main contribution is to study how a mixing condition on the covariate process transfers to a mixing condition for the response.
Using a coupling method, we first derive mixing conditions for some Markov chains in random environments, which gives a first result for some 
autoregressive categorical processes with strictly exogenous regressors. Our result is then extended to some infinite memory categorical processes. In the second part of the paper, we study autoregressive models for which the covariates are sequentially exogenous. Using a general random mapping approach on finite sets, we get explicit mixing conditions that can be checked for many categorical time series found in the literature, including multinomial autoregressive processes, ordinal time series and dynamic multiple choice models. 
We also study some autoregressive count time series using a somewhat different contraction argument.
Our contribution fill an important gap for such models, presented here under a more general form, since such a strong mixing condition is often assumed in some recent works but no general approach is available to check it.  
\end{abstract}
\vspace*{1.0cm}

\footnoterule

\section{Introduction}
Discrete-valued time series are often encountered in many real-world problems. See for instance \citet{Weiss} for an interesting textbook presenting various models and applications concerning those time series.
We usually distinguish time series taking values in a finite set $E$ which are called categorical and count time series which takes valued in the infinite set $E=\N$. In this paper, we will focus on categorical time series based on regression theory. See \citet{Fokianos2003} for a survey of the models used in this case. But we will also consider autoregressive count time series models called INGARCH models that are widely used by practitioners. 
When it comes to study semi-parametric or non-parametric estimation procedures in time series analysis, it is often necessary to derive dependence properties such as strong mixing properties from which various limit theorems can be used to derive their asymptotic properties. See for instance \citet{Doukhan(1994)},  a standard reference for this topic. 
Many models for discrete-valued time series are based on Markov chains and there exist numerous criteria to exhibit mixing properties in this case.
For instance, it is widely known that irreducible finite-state Markov chains are automatically $\phi-$mixing which entails automatically such mixing properties for the logistic/probit autoregressive models discussed in \citet{Fokianos2003}. For INGARCH processes, \citet{Neumann2011} studied $\beta-$mixing properties 
when the random intensity forms an autoregressive process with a contracting link function. \citet{doukhan2012weak} also studied some weak dependence properties for time series of counts as well as some (constrained) strong mixing coefficients.

A more tricky problem concerns strong mixing properties in such models when exogenous random covariates are included in the dynamic, which is an important problem for a realistic modeling, since in practice these models are always used with external regressors. However, the development of the theory of discrete-valued time series with exogenous covariates is quite recent and rudimentary. \citet{Truquet} considered finite-state Markov chains with strictly exogenous covariates 
and used the framework of Markov chains in random environments for deriving ergodic properties of logistic type autoregressive models.
An extension to models with infinite dependence was studied in \citet{Truquet2}, as well as some weak dependence properties.
\citet{deJ} studied a dynamic binary choice model with not necessarily strictly exogenous covariates and derived strong mixing properties of their model.
\citet{DT} studied a similar model as well as others types of autoregressive models with the functional dependence measure.
However, a systematic study of strong mixing properties for all the models listed above is still missing.

On the other hand, strong mixing properties of categorical time series are often assumed in contributions devoted to statistical inference. 
For instance, \citet{Fern} impose a strong mixing condition for studying dynamic binary/count panel data with heterogeneity.
\citet{park} considered non-parametric estimation for dynamic discrete choice models and assumed a strong mixing assumption both for the categorical response and the regressors for getting asymptotic properties. Since it is natural to include lag values of the response in the set of regressors, it could be interesting to get some examples of regression functions for which this mixing condition is indeed valid. \citet{linton} studied semi parametric inference in dynamic multiple choices models and also assumed a strong mixing conditions on both the choice and the covariates process. \citet{deJ} studied a dynamic binary choice model and used a strong mixing condition for deriving asymptotic distribution of Horowitz's smoothed maximum score estimator. 

The aim of this paper is to discuss this important issue. In particular, for categorical time series, we will show that such strong mixing conditions are often valid under quite natural 
assumptions on the model, provided that the exogenous regressors satisfy a similar condition. However such a problem cannot be studied using the classical approaches for deriving strong mixing conditions. The models are typically non-linear and do not satisfy a Markov property in general. Let us precise that 
we do not want to impose a specific dynamic structure on the covariate process. This is natural since the probability distribution of the exogenous regressors  do not afford any information for studying the dynamic of the outcome.  

 To present the setup, assume that the model can be written under the form
\begin{equation}\label{main}
Y_t=f\left(Y_{t-1},\ldots,Y_{t-p},X_{t-1},\varepsilon_t\right),
\end{equation}
where $X_{t-1}\in \R^d$ denote the set of regressors observed at previous time, $\varepsilon_t$ is a noise component and $f:E^p\times \R^d\times \R^k\rightarrow E$ is a measurable function. For ordinal time series or multiple choice models, as discussed in Section \ref{example2}, this representation is natural. On can replace $X_{t-1}$ with $(X_{t-1},\ldots,X_{t-p})$ but the latter vector can still be denoted by $Z_{t-1}$ and one can rewrite (\ref{main}) accordingly.
When the model is defined from a conditional distribution
\begin{equation}\label{main2}
\P\left(Y_t=y\vert X_{t-1},Y_{t-1},\ldots\right)=Q\left(y\vert X_{t-1},Y_{t-1},\ldots\right),
\end{equation}
one can easily go back to the representation (\ref{main}) by taking a uniformly distributed noise component $\varepsilon_t$.

However, checking our conditions require some independence assumptions between the covariates and the noise component. The most stringent one, called strict exogeneity condition, which is formally defined by the conditional independence condition (\ref{strictstrict}), is satisfied when the two processes $(X_t)_{t\in\Z}$ and $(\varepsilon_t)_{t\in\Z}$ are assumed to be independent, the latter being i.i.d.
In this case, we will use the framework of Markov chains in random environments, since conditionally on $(X_t)_{t\in\Z}$, the dynamic of $(Y_t)_{t\in\Z}$ is defined by a time-inhomogeneous Markov chains with covariates-dependent transition matrices. 
Such a framework allows to deal with a wide class of dynamics. 
But more realistic models can be considered by using a condition called predetermindness or sequential exogeneity in the econometric literature.
This condition means that $\varepsilon_t$ is independent from $\left(X_{t-j},\varepsilon_{t-j}\right)_{j\geq 1}$, with $X_t$ possibly dependent of $\varepsilon_t$. This weaker exogeneity condition is used for instance in 
\citet{PARX} for count time series and in \citet{Hsiao} (see Chapter $7.5$) for binary models for panel data.

Let us explain the general idea to get strong mixing conditions for models (\ref{main}). We first recall the definition of the strong mixing coefficients of a stationary process $(V_t)_{t\in\Z}$. For two sigma-algebras $\mathcal{F}$ and $\mathcal{G}$ on the same measurable space $\Omega$, we recall that their strong mixing coefficient is defined by
$$\alpha\left(\mathcal{F},\mathcal{G}\right)=\sup\left\{\left\vert \P(A\cap B)-\P(A)\P(B)\right\vert: (A,B)\in \mathcal{F}\times \mathcal{G}\right\}.$$
Now if $(V_t)_{t\in\Z}$ is an arbitrary stationary process taking values in taking values in an arbitrary measurable space $(H,\mathcal{H})$, we set 
$$\alpha_V(n)=\alpha\left(\mathcal{F}_V(0),\mathcal{G}_V(n)\right),$$
where for $k\in\Z$,
$$\mathcal{F}_V(k)=\sigma\left(V_j: j\leq k\right),\quad \mathcal{G}_V(k)=\sigma\left(V_j: j\geq k\right).$$
The process $(V_t)_{t\in\Z}$ is said to be strongly mixing if $\lim_{n\rightarrow \infty}\alpha_V(n)=0$. For the dynamic (\ref{main}), $V_t=(Y_t,X_t,\varepsilon_t)$. 
For a stationary dynamic of type (\ref{main}), $Y_t$ will depend on infinitely many lag values of $(\varepsilon_t,X_t)$. The main idea 
is to introduce a coupling $Y'_t$ of $Y_t$ but for which the dynamic (\ref{main}) is initialized at time $t=r_n\nearrow \infty$ where $0<r_n<n$. Setting $V_t'=(Y_t',X_t,\varepsilon_t)$, it is straightforward to show that  
\begin{equation}\label{mainbound}
\alpha_V(n)\leq \alpha\left(\mathcal{F}_V(0),\mathcal{G}_{V'}(r_n)\right)+2\sum_{t\geq n}\P\left(Y_t\neq Y'_t\right).
\end{equation}
Since the first term in the previous bound can be bounded by the mixing coefficient of the process $(X_t,\varepsilon_t)$, it remains to control the probabilities $\P\left(Y_t\neq Y'_t\right)$ which are expected to decrease to $0$ as $n\rightarrow \infty$. For categorical processes, such a control requires some care due to the non-linearity of the mapping $f$ which can be obtained from a discretization of a continuous response. See Section \ref{4} for details.
However for model of type (\ref{main2}), a good coupling will be not necessarily given by $(\ref{main})$ with a uniformly distributed noise component $\varepsilon_t$. This is why in the case of strictly exogenous regressors i.e.
\begin{equation}\label{strictstrict}
\P\left(Y_t=y\vert X_{t-1},Y_{t-1},\ldots\right)=\P\left(Y_t=y\vert (X_s)_{s\in\Z},Y_{t-1},Y_{t-2},\ldots\right),
\end{equation}
we will study other types of couplings constructed directly from the transition kernel $Q$ that defines the dynamic. Let us mention that (\ref{strictstrict}) is an exogeneity notion that can be found in \citet{Chamb}.
The reason for separating the study of strict and sequential exogeneity is mainly technical. Though more restrictive, strict exogeneity allows to work conditionally on a covariate process which plays the rule of a random environment. One can then use quite powerful coupling techniques with a very general model specification. On the other hand, sequential exogeneity require to control specific stochastic iterations that are dependent across the time, due to the presence of the covariate process. In this case, the techniques used are different from that only involving strictly exogenous regressors. In the present paper, this distinction will be only made for finite state spaces. For autoregressive count time series, \cite{DNT} recently introduced very general results for getting existence of a stationary and ergodic count time series when strictly exogenous regressors are incorporated in the dynamic. But the unbounded state space is more difficult to tackle and we did not find a way for controlling the mixing coefficients without stronger assumptions on the model. Let us mention that throughtout this paper, we focus on stationary models. Our result could be extended to accomodate with non-stationary covariates, provided that all the quantities required 
in our bounds can be made uniform with respect to the time $t$. However, this framework will not cover many interesting non-stationary processes, such as the the non-stationary binary choice model of \citet{Ph}, based on a unit-root covariate process. We then prefer to restrict to the stationary case, for conciseness of the exposure.

The paper is organized as follows. In Section \ref{2}, we consider Markov chains models in random environments which is adapted to study dynamics of type (\ref{main2}) depending of finitely many past values under a strict exogeneity assumption. The dependence with respect to infinitely many past values is investigated in Section \ref{3} and covers some observation-driven models found in the literature.
In Section \ref{4}, we will switch to sequential exogeneity in (\ref{main}) using a random mapping framework. Section \ref{5} will be devoted to infinite dependence with much more restrictive conditions on the model than in the case of strict exogeneity. The case of binary models and INGARCH type processes are then investigated. Finally an Appendix Section \ref{6} provides some useful lemmas for our proofs.

\section{Mixing properties for Markov chains in random environments}\label{2}
In this section, we establish strong mixing properties for the models investigated in \citet{Truquet}. Though our main interest is 
to study finite-state Markov chains $(Y_t)_{t\in\Z}$ defined conditionally on a strictly exogenous covariate process $(X_t)_{t\in\Z}$, 
we give a more general result for Markov chain models in random environments satisfying some Doeblin's type condition. 
Such models have been considered by \citet{Kifer} with a slightly more general structure. In what follows, let $E$ and $F$ be two Polish spaces and
$\left\{P_x: x\in F\right\}$ be a family of Markov kernels on a Polish space $E$. 
We assume that the mappings $(x,y)\mapsto P_x(y,A)$ are measurable for any $A\in\mathcal{B}(E)$.
In what follows, we consider a stochastic process $(X_t)_{t\in \Z}$, called the environment, taking values in $F$ and
$(Y_t)_{t\geq 0}$ a sequence of $E-$valued random variables such that a.s.
\begin{equation}\label{main1}
\P\left(Y_t\in A\vert X,Y_{t-1},\ldots,Y_0\right)=P_{X_{t-1}}\left(Y_{t-1},A\right).
\end{equation}
The process $(Y_t)_{t\geq 0}$ is called Markov chain in random environments. Under some conditions, there exists a unique stationary solution $\left((X_t,Y_t)\right)_{t\in\Z}$
for (\ref{main1}). In what follows, we fix a positive integer $m$ and we set $Z_t=\left(X_t,X_{t+1},\ldots,X_{t+m-1}\right)$ for any $t\in\Z$.

\begin{description}
\item[A1] The process $(X_t)_{t\in\Z}$ is stationary and strongly mixing.
\item[A2] The following Doeblin's type condition is satisfied. We have almost surely,
$$P_{X_1}\cdots P_{X_m}(x,y)\geq \eta_{Z_1}\nu_{Z_1}(y),\quad (x,y)\in E^2,$$
where $z\mapsto \eta_z$ is a measurable mapping from $F^m$ to $(0,1)$ and for each $z\in F^m$, $\nu_z$ is a probability measure on $E$.
\end{description}

Our aim is to derive mixing conditions for the process $\left((X_t,Y_t)\right)_{t\in \N}$ when the process $(X_t)_{t\in \N}$ is itself mixing.
When Assumption {\bf A2} holds true and $(X_t)_{t\in\Z}$ is stationary, \citet{Kifer} showed the existence of a unique stationary distribution for such a problem. In particular, for any $y\in E$,
the sequence of random measures $\left(\nu_n^{(y)}\right)_{n\geq 1}$ defined by $\nu_n^{(y)}=P_{X_{t-n}}\cdots P_{X_{t-1}}\left(y,\cdot\right)$ converges almost surely to a random measure $\pi_t$ in the total variation sense.
Moreover $\pi_t$ does not depend on $y$. It is straightforward to show that 
$$\P\left(Y_t\in A\vert (X_t)_{t\in\Z}\right)=\pi_t(A).$$
Moreover, the pair $(X_t,Y_t)$ is ergodic provided that the environment is itself ergodic. Let us also mention that \cite{HL} considered convergence of similar backward products $\nu_n^{(y)}$ under some weaker assumptions  but {\bf A2} is sufficient for the examples we will discuss.
The main result of this section is given below.

\begin{theo}\label{mixing1}
Suppose that Assumptions {\bf A1-A2} hold true. Set $\rho=1-\E\eta_{Z_0}$.
\begin{enumerate}
\item
Then the strong mixing coefficients of the pair $V_t=\left(X_t,Y_t\right)$ are bounded as follows. For any pair of integers $r\leq t$ such that $1\leq r\leq n-1$, set $s_t(r)=[(t-r)/m]$. We then have for any $1\leq r\leq n-1$,
$$\alpha_V(n)\leq 4 \alpha_X(r)+2\sum_{t\geq n}\inf_{1\leq j\leq s_t(r)-1}\left\{\rho^{[s_t(r)/j]}+4\frac{\alpha_X\left((j-1)m+1\right)}{1-\rho}\right\}.$$
\item
In particular, assume that $\alpha_X(n)=O\left(n^{-\kappa}\right)$ with $\kappa>1$. Then, for any real number $\kappa'$ such that $1<\kappa'<\kappa$, the process $(V_t)_{t\in\Z}$ is strongly mixing with $\alpha_V(n)=O\left(n^{-\kappa'+1}\right)$.
If now $\alpha_X(n)=O\left(\kappa^n\right)$ with $\kappa\in (0,1)$, there then exists $\overline{\kappa}\in (0,1)$ such that $\alpha_V(n)=O\left(\overline{\kappa}^{\sqrt{n}}\right)$. 
\end{enumerate}
\end{theo} 

\paragraph{Notes}
\begin{enumerate}
\item
When $E$ is finite, say $E=\{1,\ldots,N\}$, and $F=\R^d$, there is a generic way for constructing random stochastic matrices satisfying {\bf A2}. Take a regular transition matrix $P$ (i.e. there exists a power $m$ of $P$ with positive entries).
If $J_i$ is the subset of $E$ such that $P(i,j)>0$ for $j\in J_i$, define 
$$P_{X_t}(i,j)=\exp\left(\theta_{i,j}'X_t\right)/\left(\sum_{\ell\in J_i}\exp\left(\theta_{i,\ell}'X_t\right)\right),$$
where $\theta_{i,\ell}$ is a column vector of $\R^d$ and $\theta_{i,\ell}'$ denotes its transpose. One can check that Assumptions {\bf A2} is automatically satisfied.
\item
If $\eta:=\eta_{Z_1}$ is deterministic, i.e. the constant in the Doeblin's condition {\bf A2} is uniform with respect to the environment, inspection of our proof shows that the quantity
$$\inf_{1\leq j\leq s_t(r)-1}\left\{\rho^{[s_t(r)/j]}+4\frac{\alpha_X\left((j-1)m+1\right)}{1-\rho}\right\}$$
can be replaced with $\left(1-\eta\right)^{s_t(r)}$. Hence $\alpha_V(n)=O\left(\alpha_X(n)\right)$ if $\alpha_X(n)$ has a power decay. Note also that $\alpha_V(n)$ has a geometric decay as soon as the same property holds true for $\alpha_X(n)$. For most interesting models, a deterministic constant can only be obtained when the environment forms a bounded process.
\item
In the finite-state case, our result applies to some models defined by 
$$\P\left(Y_t=y\vert (X_t)_{t\in\Z}, Y_{t-1},\ldots\right)=H_y\left(X_{t-1},Y_{t-1},\ldots,Y_{t-p}\right),$$
with $\left\{H_y: y\in E\right\}$ a family of measurable functions, taking values in $(0,1)$ and such that $\sum_{y\in E}H_y=1$. Indeed, {\bf A2} is satisfied with 
$m=p$ and 
$$P_{X_t}\left((y_1,\ldots,y_p),(y_2,\ldots,y_{p+1})\right)=H_{y_{p+1}}\left(X_t,y_p,\ldots,y_1\right),$$
which corresponds to the transition matrix of $\left(Y_t,Y_{t+1},\ldots,Y_{t+p-1}\right)$. We will revisit this model without the strict exogeneity condition. See (\ref{mainex2}) for details.

\item
Our results can be also useful for other dynamics than finite-state processes. For instance, consider the case $E=[0,1]$, $F=\R^d$ and a beta autoregressive process $(Y_t)_{t\in \Z}$ with
$P_x(y,\cdot)$ being the beta distribution with positive parameters $\left(a_1(x,y),a_2(x,y)\right)$ and such that for $i=1,2$,
$$0<\inf_{y\in E}a_i(x,y)\leq \sup_{y\in E} a_i(x,y)<\infty.$$
One can show that {\bf A2} is satisfied with $m=1$ and $\nu_{X_1}$ being the beta distribution with parameters $\left(\sup_{y\in E}a_1(X_1,y),\sup_{y\in E}a_2(X_1,y)\right)$ and $\eta_{X_1}$ a suitable positive random variable. See \citet{beta} for some models of this type.

\end{enumerate}

\paragraph{Proof of Theorem \ref{mixing1}}
\begin{enumerate}
\item
Our aim is to derive a bound of type (\ref{mainbound}) for a suitable version $Y'$ of the Markov chain in random environments.
Let $1<r<n$. We set $Y_j'=y_0$ where $y_0$ is an arbitrary point in $E$. On a triplet $\left(\Omega,\mathcal{A},\P\right)$ on which the covariate process $(X_t)_{t\in\Z}$ is defined, we assume that the pair $\left(Y_t,Y'_t\right)_{t\in \Z}$ is defined as follows. 
First $(Y_t)_{t\leq r}$ is such that $\left(X_t,Y_t\right)_{t\leq r}$ is a stationary process such that (\ref{main1}) holds true.
Let $y_0$ be an arbitrary point in $E$ and we denote by $\mathcal{G}_r$ the sigma-field generated by $(X_t)_{t\in\Z}$ and $(Y_t)_{t\leq r}$.
We then assume that $\P\left(Y_r'=y_0\vert \mathcal{G}_r\right)=1$. Now for any $t\in \Z$, we denote by $R_{Z_t}$ the random probability kernel defined by
$$R_{Z_t}(y,A)=\frac{P_{X_t}\cdots P_{X_{t+m-1}}(y,A)-\eta_{Z_t}\nu_{Z_t}(A)}{1-\eta_{Z_t}},\quad (y,A)\in E\times \mathcal{B}(E)$$
and next for $y,y'\in E$ and $A,B\in \mathcal{B}(E)$,
$$Q_{Z_t}\left(y,y';A\times B\right)\eta_{Z_r}\nu_{Z_r}(A\cap B)+\left(1-\eta_{Z_r}\right)R_{Z_r}\left(y,A\right)R_{Z_r}\left(y',B\right)$$
We then assume that for $(A,B)\in \mathcal{B}(E)^2$,
\begin{eqnarray*}
\P\left(Y_{m+r}\in A,Y'_{m+r}\in B\vert \mathcal{G}_r\right)&=&\mathds{1}_{\{Y_r\neq Y_r'\}}Q_{Z,r}\left(Y_r,Y'_r; A\times B\right)\\&+&\mathds{1}_{\{Y_r=Y'_r\}}P_{X_r}\cdots P_{X_{r+m-1}}\left(Y_r,A\cap B\right).
\end{eqnarray*}
The interpretation of the coupling is as follows. Conditionally on past values, if $Y_r\neq Y_r'$, with probability $\eta_{Z_r}$, we draw $Y_{m+r}=Y'_{m+r}\sim \eta_{Z_r}$ and with probability $1-\eta_{Z_r}$, we draw two independent random variables 
with the transition kernel $R_{Z_r}$. If $Y_r=Y'_r$, we draw $Y_{m+r}=Y'_{m+r}\sim P_{Z_r}\left(Y_r,\cdot\right)$. 
Next, conditionally on $\mathcal{G}_r\vee \sigma\left(Y_{r+m},Y'_{r+m}\right)$, we simulate $\left(Y_{r+1},\ldots,Y_{r+m-1}\right)$ and $\left(Y'_{r+1},\ldots,Y'_{r+m-1}\right)$ independently according to the distribution 
$W_{Z_r}\left(Y_r,Y_{r+m}; \cdot\right)$ and $W_{Z_r}\left(Y'_r,Y'_{r+m}; \cdot\right)$, where $W_{z}\left(y_0,y_m;\cdot\right)$ is a version of the conditional distribution of $\left(Y_{r+1},\ldots,Y_{r+m-1}\right)$ given that 
$Z_r=z,Y_r=y_0,Y_{r+m}=y_m$. 
We next draw a new block of $m$ successive values $\left(Y_{m+r+\ell},Y'_{m+r+\ell}\right)_{1\leq \ell\leq m}$ using the same strategy and so on.
It is easy to check that $(Y_t)_{t\in \Z}$ and $(Y'_t)_{t\geq r}$ are two paths of the same Markov chain in random environments, the latter one being initialized at time $t=r$ with $Y_r'=y_0$.  
Now if $t=r+j+sm$ with $0\leq j\leq m-1$ and $s\in \N^{*}$, we have
\begin{eqnarray*}
\P\left(Y_t\neq Y'_t\right)&\leq &\P\left(Y_{r+s m}\neq Y'_{1+s m},\ldots,Y_{r+m}\neq Y'_{r+m}\right)\\
&\leq& \E\left(\prod_{\ell=0}^{s-1}\left(1-\eta_{Z_{r+\ell m}}\right)\right).
\end{eqnarray*}
Using Lemma \ref{ult} and remembering that $\rho=1-\E\eta_{Z_0}$, we get 
$$\P\left(Y_t\neq Y'_t\right)\leq \inf_{1\leq r\leq s-1}\left\{\rho^{[s/r]}+\frac{\alpha_X\left((r-1)m-1\right)}{1-\rho}\right\}.$$
Note that 
$$\alpha_V(n)\leq \alpha\left(\mathcal{F}_0,\mathcal{G}'_n\right)+2\sum_{t\geq n}\P\left(Y_t\neq Y'_t\right).$$
Set $X:=(X_t)_{t\in\Z}$. We have for two events $A$ and $B$ in the cylinders sigma-field,
$$\P\left((Y_t,X_t)_{t\leq 0}\in A, (Y'_t,X_t)_{t\geq n}\in B \vert X\right)=\P\left((Y_t,X_t)_{t\leq 0}\in A\vert X\right)\times \P\left((Y_t,X_t)_{t\leq n}\in B\vert X\right).$$
Indeed $(Y_t')_{t\geq r}$ has been constructed independently from $(Y_t)_{t\leq 0}$, conditionally on $X$.
Moreover, $G_A:=\P\left((Y_t,X_t)_{t\leq 0}\in A\vert X\right)$ is a measurable function of $(X_t)_{t\leq 0}$ and $H_B:=\P\left((Y'_t,X_t)_{t\geq n}\in B\vert X\right)$ is a measurable function of $(X_t)_{t\geq r}$.
We then have 
$$\alpha\left(\mathcal{F}_0,\mathcal{G}'_n\right)\leq \sup_{A,B}\left\vert\c\left(G_A,H_B\right)\right\vert\leq 4 \alpha_X(r),$$
where the last inequality follows from the covariance inequality given in \citet{Doukhan(1994)}, Lemma $3$. The announced upper-bound is now proved.
\item
In the case of power decays of the mixing coefficients, we choose $r\sim n/2$ and then $j\sim s_t(r)^{\ell}$ for some $\ell\in (0,1)$ such that $\kappa'=\ell \kappa>1$ and apply the previous point. 
For a geometric decay of the mixing coefficients of $\zeta$, we also set $r\sim n/2$ and $j=\sqrt{s_t(r)}$ and the result can obtained by noticing that for any $\overline{\rho}\in (0,1)$,
$\sum_{t\geq k}\overline{\rho}^{\sqrt{t}}=O\left(\widetilde{\rho}^{\sqrt{t}}\right)$ for any $\widetilde{\rho}>\overline{\rho}$, if we use a comparison between the series and an integral.$\square$
\end{enumerate}

\section{Mixing properties for infinite memory models with strictly exogenous regressors}\label{3}
In this section, we extend the result of the previous one, when the number of regressors is infinite. 
In \citet{Truquet2}, such models has been considered and weak dependence properties such as $\beta-$mixing properties have been derived. 
The result presented below is more general and applies to $\alpha-$mixing covariates. Moreover, we will not assume that the covariate process
is a function of a Markov chain. 
The framework used in this part is based on the theory of chains with complete connections. See \citet{Chaz} and the reference therein for a recent contribution to these models and the battery of results developed to study them. 
We consider models defined by 
\begin{equation}\label{maininfinite}
\P\left(Y_t=y\vert Y_{t-1}^{-},X_{t-1}^{-}\right)=P\left(y\vert Y_{t-1}^{-},X_{t-1}^{-}\right),
\end{equation}
where $P$ is a transition kernel from $E^{\N}\times \mathcal{D}$ to $E$ where $\mathcal{D}$ is a measurable set of $(\R^d)^{\N}$ such that $\P\left(X_{-1}^{-}\in \mathcal{D}\right)=1$.
Here, for any sequence $x\in (\R^d)^{\Z}$, $x_t^{-}$ denotes the sequence $\left(x_{t-j}\right)_{j\geq 0}$. 
As explained in \citet{Truquet2}, models of type (\ref{maininfinite}) include many observation-driven categorical times series models found in the literature, in particular in Econometrics. See for instance \citet{kauppi} or \citet{rydberg} or \citet{russell}.
In what follows, we denote by $d_{TV}(\mu,\nu)$ the total variation distance between two probability measures $\mu$ and $\nu$ on $E$, i.e.
\begin{eqnarray*}
d_{TV}(\mu,\nu)&=&\frac{1}{2}\sum_{y\in E}\left\vert \mu(y)-\nu(y)\right\vert\\
&=&1-\sum_{y\in E}\min\left(\mu(y),\nu(y)\right).
\end{eqnarray*}
We also consider an arbitrary norm $\vert\cdot\vert$ on $\R^d$ and set $\Vert X_0\Vert_1=\E\vert X_0\vert$.
We will use the two following assumptions on the transition kernel.

\begin{description}\label{supersuper}
\item[A3] We have
$$\sup_{y,y'\in E^{\N}}\sup_{x\in \mathcal{D}}d_{TV}\left(P\left(\cdot\vert y,x\right),P\left(\cdot\vert y',x\right)\right)<1.$$
\item[A4] There exist two sequences of nonnegative real numbers $(a_j)_{j\geq 0}$ and $(e_0)_{j\geq 1}$ such that $\sum_{j\geq 0}j a_j<\infty$, $\sum_{j\geq 0}e_j<\infty$ and for $(y,y',x,x')\in E^{\N}\times E^{\N}\times \mathcal{D}\times\mathcal{D}$, 
$$d_{TV}\left(P\left(\cdot\vert y,x\right),P\left(\cdot,y',x'\right)\right)\leq \sum_{i\geq 0}a_i \mathds{1}_{y_i\neq y_i'}+\sum_{j\geq 0}e_j\vert x_j-x_j'\vert.$$
\end{description}
These assumptions guarantee the validity of Assumptions {\bf S2-S3} in \citet{Truquet2}. Indeed, the coefficients used in this reference are given by 
$$b_m=\sup_{y,y'\in E^{\N},x\in \mathcal{D}}\left\{d_{TV}\left(P\left(\cdot\vert y,x\right),P\left(\cdot\vert y',x\right)\right): y_j=y'_j, 0\leq j\leq m-1\right\},$$
defined for a positive integer $m$ satisfies the bound $b_m\leq \sum_{j\geq m+1}a_j$ and are then summable.
When $(X_t)_{t\in\Z}$ is a stationary and ergodic process, Theorem $1$ in \citet{Truquet2} entails the existence and uniqueness (in the probability distribution sense) of a stationary and ergodic process $\left((Y_t,X_t)\right)_{t\in\Z}$ solution of (\ref{maininfinite}).
When the kernel $P$ depends of finitely many values, we recover the setup of the Markov chains in random environments of the previous section.
But the assumptions used here are more restrictive. Takes for instance the case $m=1$ and assume that 
$P\left(\cdot\vert Y_{t-1}^{-},X_{t-1}^{-}\right)=P_{X_{t-1}}\left(Y_{t-1},\cdot\right)$. Assumption {\bf A3} is based on a control of the total variation distance which is uniform with respect to the covariates, while the random Doeblin's condition is {\bf A2} is compatible with a non-uniform control of such total variation distances. Indeed, under {\bf A2},
$$d_{TV}\left(P_{X_{t-1}}\left(y,\cdot\right),P_{X_{t-1}}\left(y,\cdot\right)\right)\leq 1-\eta_{X_{t-1}}.$$  
However, we now cover the case of infinite dependence which is substantially more difficult than the Markov case.
In what follows, we set $S_j=\sum_{t\geq j}e_j$ for any nonnegative integer $j$ and for two sequences of real numbers $(u_j)_{j\geq 0}$ and $(v_j)_{j\geq 0}$, we denote by $u*v$ their convolution product, i.e. $(u*v)_n=\sum_{j=0}^n u_j v_{n-j}$.
We get the following result.

\begin{theo}\label{mixing2}
Suppose that {\bf A3-A4} hold true and that $(X_t)_{t\in\Z}$ is stationary and ergodic. 
For any integer $r$ such that $1<r<n$, we then have the bound
$$\alpha_V(n)\leq 4\alpha_X(r)+2\sum_{t\geq n-r-1}b_t^{*}+4\Vert X_0\Vert_1\sum_{t\geq n-r}S_t+4\Vert X_0\Vert_1\sum_{t\geq n-r-1}\left(b^{*}*S\right)_t.$$
The sequence $(b^{*}_j)_{j\geq 0}$ is defined by $b_0^{*}=b_0$ and for $n\geq 1$, $b_n^{*}$ is equal to $\P\left(T_n^{(b)}=0\right)$ where $\left(T_n^{(b)}\right)_{n\geq 0}$ is a time-homogeneous Markov chain, starting at $0$ and with transition matrix $Q$ defined by 
$$Q(i,i+1)=1-b_i,\quad Q(i,0)=b_i,\quad i\in\N.$$
\end{theo}
 
\paragraph{Note.} Take for instance the model with $E=\{0,1\}$ and
$$\P\left(Y_t=1\vert X_{t-1}^{-},Y_{t-1}^{-}\right)=F\left(\lambda_t\right),\quad \lambda_t=\beta \lambda_{t-1}+\kappa Y_{t-1}+\delta' X_{t-1}.$$
Here $F$ denotes a cumulative distribution function with full support such that the Gaussian c.d.f. (probit model) or the logistic (i.e. $F(s)=(1+e^{-s})^{-1}$).
From Proposition $3$ in \citet{Truquet2}, existence and uniqueness of a stationary and ergodic solution is guaranteed as soon as $\vert a\vert <1$ and $\E\vert X_0\vert<\infty$.  In this case, the coefficients $b_m$ decays geometrically and so do the coefficients $b_m^{*}$ and $\gamma_{\ell}$. Taking $r=[n/2]$ in the previous results, we get 
$$\alpha_V(n)\leq C\left\{\alpha_X([n/2])+ \rho^n\right\},$$
for some constants $(C,\rho)\in (0,\infty)\times (0,1)$. The same kind of bound can be obtained for more general observation-driven models, as discussed in \citet{Truquet2}.

\paragraph{Proof of Theorem \ref{mixing2}}
Let $1<r<n$ and $y_0$ be an arbitrary state in $E$. Set $Y'_t=y_0$ for $t\leq r$. Suppose that $(X_t)_{t\in\Z}$ and $(Y_t)_{t\leq r}$ are already defined with $(Y_t,X_t)_{t\leq r}$ stationary as well as $\left(\widetilde{X}_t\right)_{t\leq r-1}$ a copy of $(X_t)_{t\leq r-1}$, independent of $\sigma\left((X_j,Y_s): j\in \Z, s\leq r\right)$. We set $\widetilde{X}_t=X_t$ for $t\geq r$. From Lemma $1$ in \citet{Truquet2}, it is possible to construct a $(Y_t,Y'_t)_{t\geq r+1}$ such that 
$$\P\left(Y_t=y\vert Y_{t-1}^{-},X_{t-1}^{-}\right)=P\left(\cdot\vert Y_{t-1}^{-},X_{t-1}^{-}\right),\quad P\left({Y'}_t=y\vert {Y'}_{t-1}^{-},\widetilde{X}_{t-1}^{-}\right)=P\left(\cdot\vert {Y'}_{t-1}^{-},X_{t-1}^{-}\right),$$
for all $y\in E$ and $t\geq r+1$ and
\begin{eqnarray*}
\P\left(Y_t\neq {Y'}_t\vert X\right)&\leq& b^{*}_{t-r-1}+\sup_{s\in E^{\N}}d_{TV}\left(P\left(\cdot\vert s,X_{t-1}^{-}\right),P\left(\cdot\vert s, \widetilde{X}_{t-1}\right)\right)\\
&+& \sum_{\ell=0}^{t-r-2} b_{\ell}^{*} \sup_{s\in E^{\N}}d_{TV}\left(P\left(\cdot\vert s, X_{t-\ell-2}\right),P\left(\cdots,s, \widetilde{X}_{t-\ell-2}^{-}\right)\right).
\end{eqnarray*}
Using {\bf A4}, we conclude after a few computations that
$$\P\left(Y_t\neq Y'_t\right)\leq b_{t-r-1}^{*}+2\E\vert X_0\vert\times\left(S_{t-r}+ \sum_{\ell=0}^{t-r-2}b_{\ell}^{*}S_{t-r-\ell-1}\right).$$
As in the proof of Theorem \ref{mixing1}, we obtain the following bound
$$\alpha_V(n)\leq \alpha\left(\mathcal{F}_0,\mathcal{G}'_r\right)+2\sum_{t\geq n}\P\left(Y_t\neq Y'_t\right).$$
The end of the proof is similar to that of Theorem \ref{mixing2}. Indeed, from Theorem $1$ in \citet{Truquet2}, the conditional distribution of $Y_0^{-}$ given $(X_t)_{t\in\Z}$ only depend on $X_0^{-1}$. 
In what follows, we set $x^{+}_t=(x_{t+j})_{j\geq 0}$ for any sequence $(x_t)_{t\in \Z}$ of elements in an arbitrary space.
We then get for two events $A$ and $B$ in the cylinder sigma-field,
$$\P\left(V_0^{-}\in A\vert (X_t)_{t\in \Z},(\widetilde{X}_t)_{t\in\Z}\right)=\P\left(V_0^{-}\in A\vert X_0^{-}\right),$$
$$\P\left({V'}_r^{+}\in B\vert (X_t)_{t\in \Z},(\widetilde{X}_t)_{t\in\Z}\right)=\P\left({V'}_r^{+}\in B\vert X_r^{+},\widetilde{X}^{-}_{r-1}\right).$$
Using the independence between $(X_t)_{t\in\Z}$ and $\widetilde{X}^{-}_{r-1}$ and the covariance inequality for strong mixing variables, we get 
$$\alpha\left(\mathcal{F}_0,\mathcal{G}'_r\right)\leq 4\alpha_X(r)$$
and the proposed bound easily follows.$\square$

\section{Mixing properties of iterated random maps on a finite state space}\label{4}

\subsection{A general result for iterated dependent random maps}
Let $(\zeta_t)_{t\in\Z}$ be a stationary process taking values in $\R^e$. For any $s\in E$, we consider a mapping $F_s:E\rightarrow E$ and we assume 
that the mapping $(x,s)\mapsto F_s(x)$ is measurable as an application from $E\times \R^e$ to $E$.
We next consider a stochastic process $(Y_t)_{t\geq 0}$, taking values in $E$, and such that 
$$Y_t=F_{\zeta_t}(Y_{t-1}),\quad t\geq 1.$$
For simplicity of notations, we set $F_s^t=F_{\zeta_t}\circ F_{\zeta_{t-1}}\circ\cdots\circ F_{\zeta_s}$ for $s<t$. 
The following assumption will be crucial. In what follows, we denote by $\# A$ the cardinality of a set $A$.

\begin{description}
\item[B1] There exists a positive integer $m$ such that $1-\rho:=\P\left(\# F_1^m(E)=1\right)>0$.
\end{description}

Assumption {\bf B2} means that if the random maps are iterated sufficiently, there is a positive probability to get coalescence of the iterations and on the corresponding event, the system loses its memory with respect to the initial state. Note that the random maps are not independent here.

The main result of this section is the following.

\begin{theo}\label{mainit}
Suppose that the process $(\zeta_t)_{t\in\Z}$ is stationary and ergodic an that Assumption {\bf B1} holds true. 
\begin{enumerate}
\item
There then exists a unique stationary and ergodic process $(Y_t)_{t\in\Z}$ 
taking values in $E$ and such that $Y_t=F_t\left(Y_{t-1}\right)$ a.s.
\item
Moreover if the process $(\zeta_t)_{t\in\Z}$ is strongly mixing, then the strong mixing coefficients of the pair $V_t=\left(\zeta_t,Y_t\right)$ are bounded as follows. For any pair of integers $r\leq t$ such that $1\leq r\leq n-1$, set $s_t(r)=[(t-r)/m]$. We then have for any $1\leq r\leq n-1$,
$$\alpha_V(n)\leq \alpha_{\zeta}(r+1)+2\sum_{t\geq n}\inf_{1\leq j\leq s_t(r)-1}\left\{\rho^{[s_t(r)/j]}+4\frac{\alpha_{\zeta}\left((j-1)m\right)}{1-\rho}\right\}.$$
\item
In particular, assume that $\alpha_{\zeta}(n)=O\left(n^{-\kappa}\right)$ with $\kappa>1$. Then, for any real number $\kappa'$ such that $1<\kappa'<\kappa$, the process $(V_t)_{t\in\Z}$ is strongly mixing with $\alpha_V(n)=O\left(n^{-\kappa'+1}\right)$.
If now $\alpha_{\zeta}(n)=O\left(\kappa^n\right)$ with $\kappa\in (0,1)$, there then exists $\overline{\kappa}\in (0,1)$ such that $\alpha_V(n)=O\left(\overline{\kappa}^{\sqrt{n}}\right)$. 
\end{enumerate}
\end{theo}

\paragraph{Proof of Theorem \ref{mainit}}
\begin{enumerate}
\item
For the first statement, we simply show that the random sequence $\left(F_{t-n}^t(y)\right)_{n\geq 0}$ has an almost sure limit not depending on $y$ when $n\rightarrow \infty$.
Our argument is based on a coalescence argument already used for specific dynamics by \citet{DT}.
However, the approach used here is much more synthetic and it will be convenient for including a wider class of models.
Of course, since the state space is finite, convergence means that these iterations are constant when $n$ is large enough.
Set $Z_t=\left(\zeta_t,\ldots,\zeta_{t+m-1}\right)$ for $t\in\Z$.
Since the process $\left(Z_{t-j},\zeta\right)_{j\geq 0}$ is also stationary and ergodic, the events $A_{t-j-m+1}^{t-j}:=\left\{\# F_{t-j-m+1}^{t-j}(E)=1\right\}\in \sigma\left(Z_{t-j-m+1}\right)$ occur infinitely often. This is a consequence of Birkhoff's ergodic theorem. Said differently,
$$\P\left(\sum_{j=0}^{\infty}\mathds{1}_{A_{t-j-m+1}^{t-j}}=\infty\right)=1,$$  
 If $T=T_t$ defines the first (random) integer $j\geq 1$ such that $\mathds{1}_{A_{t,j}}=1$ a.s., we have 
we have for any $y_0\in E$, $F_{t-n}^t(y)=F_{t-T-m+1}^t(y_0)$ on the event $\left\{n\geq T+m-1\right\}$. This a consequence of the equalities $F_{t-T-m+1}^{t-T}(y)=F_{t-T-m+1}^{t-T}(y')$ a.s. for $y\neq y'$. This shows the announced property. 
Denoting by $Y_t$ this limit, it is quite clear that $Y_t$ writes as a measurable function $H:E^{\N}\rightarrow E$ of $\left(\zeta_{t-j}\right)_{j\geq 0}$ which shows that the process $(Y_t)_{t\in \Z}$ is stationary and ergodic. 
We next show that $Y_t=F_t(Y_{t-1})$ a.s. We have $Y_t=F_{t-T-m+1}^t(y_0)=F_t\circ F_{t-T-m+1}^{t-1}(y_0)=F_t\left(Y_{t-1}\right)$, since $Y_{t-1}=F_{t-T-m+1}^{t-1}(y_0)$ a.s.  
Now, if $\left(Y'_t\right)_{t\in \Z}$ is satisfies $Y_t'=F_t(Y_{t-1}')$ for $t\in\Z$, we have $Y_t'=F_{t-T-m+1}^t\left(Y'_{t-T-m}\right)=F_{t-T-m+1}^t(y_0)=Y_t$ a.s.    
\item
We next study the mixing properties of the unique stationary solution. 
To this end, let $0<r<n$. We set $Y'_r=y_0$ and $Y'_t=F_t\left(Y'_{t-1}\right)$ for $t\geq r+1$. We are first going to control $\P\left(Y'_t\neq Y_t\right)$.
We have for $t=j+r+sm$ with $0\leq j\leq m-1$, 
$$\left\{Y_t\neq Y'_t\right\}\subset\left\{Y_{r+sm}\neq Y'_{r+sm}\right\}=\cap_{i=1}^s\left(\Omega\setminus A_{r+(i-1)m+1}^{r+im}\right).$$
Let $\kappa_i$ be the indicator function of the event $\Omega\setminus  A_{r+(i-1)m+1}^{r+im}$ for $1\leq i\leq s$ and $\rho=\E\kappa_1$. From {\bf B1}, we have $\rho<1$. 
We then have 
$$\P\left(Y_t\neq Y'_t\right)\leq \E\left(\kappa_1\cdots \kappa_s\right).$$ 
Observe that $\alpha_{\kappa}(j)\leq \alpha_X\left((j-1)m\right)$ for $j\geq 1$. Using Lemma \ref{ult}, we then get the bound
\begin{equation}\label{supercoup}
\P\left(Y_t\neq Y'_t\right)\leq \inf_{1\leq j\leq s-1}\left\{\rho^{[s/j]}+\frac{\alpha_X\left((j-1)m\right)}{1-\rho}\right\}.
\end{equation}
We now show how such a control can be used to bound the mixing coefficients $\alpha_{\zeta,Y}$. 
Due to the Bernoulli shift representation of $Y_t$, we first note that $\alpha_V(n)\leq \alpha\left(\mathcal{F}_{\zeta}(0),\mathcal{G}_V(n)\right)$.
Now let $r$ be an integer between $1$ and $n-1$. Denote by $(Y'_t)_{t\geq r}$ the process defined by $Y'_s=y_0$ and $Y'_t=F_t\left(Y'_{t-1}\right)$ for $t\geq r+1$. From the definition of the mixing coefficients, it is quite clear that the following bound holds true.
$$\alpha_V(n)\leq \alpha\left(\mathcal{F}_{\zeta}(0),\mathcal{G}_{\zeta}(r+1)\right)+2\sum_{t\geq n}\P\left(Y_t\neq Y'_t\right).$$
The probabilities $\P\left(Y_t\neq Y'_t\right)$ can be bounded from  (\ref{supercoup}), if we take care to replace $s$ with $s_t(r)$. We then obtain the bound given in the second point of the lemma.
\item
The proof is similar to that of Theorem \ref{mixing1}, point $2$.$\square$
\end{enumerate}

\subsection{Examples}\label{example2}
In this subsection, we give many examples of autoregressive time series models for which our results can be applied. 
For simplicity, we will always denote by $G=\{1,\ldots,N\}$ the state space of the time series.

All the models considered below will write as 
\begin{equation}\label{mainex}
Y_t=f\left(Y_{t-1},\ldots,Y_{t-p},X_{t-1},\varepsilon_t\right),
\end{equation}
where $f:G^p\times \R^d\times \R^k\rightarrow G$ is a measurable function, $(X_t)_{t\in\Z}$ is a sequence of random variables taking values in $\R^d$ and $\left(\varepsilon_t\right)_{t\in\Z}$ is a noise process.
For recovering an iterated random maps system, one can set $E=G^p$, $\zeta_t=\left(X_{t-1},\varepsilon_t\right)$ and 
$$F_t\left(y_1,\ldots,y_p\right)=\left(f(y_1,\ldots,y_p,X_{t-1},\varepsilon_t),y_1,\ldots,y_{p-1}\right).$$
Obviously, $(Y_t)_{t\in\Z}$ is solution of (\ref{mainex}) if and only if $Z_t=F_t\left(Z_{t-1}\right)$, $t\in \Z$, $Z_t=\left(Y_t,\ldots,Y_{t-p+1}\right)$. Moreover we have $\alpha_Y(n)\leq \alpha_Z(n)$ and $\alpha_Z(n)$
can be bounded from Theorem \ref{mainit}. Since the mixing condition of $\zeta$ is a simple assumption to make for applying our result, we mainly concentrate 
on the checking of {\bf B1} for the different models.
Throughout this subsection, we will assume that
\begin{description}
\item[B2]
The process $\left((X_t,\varepsilon_t)\right)_{t\in\Z}$ is stationary and ergodic.
\end{description}
Note that from {\bf B2}, the process $\left(\zeta_t\right)_{t\in\Z}$ is also stationary and ergodic.
The notion called predetermindness or sequential exogeneity is stated as follows.
\begin{description}
\item[B3]
For any $t\in \Z$, $\varepsilon_t$ is independent of $\sigma\left((X_{t-j},\varepsilon_{t-j}): j\geq 1\right)$.
\end{description}
Contrarily to the case of Markov chains in random environments, $X_t$ may be dependent of $\varepsilon_t$. We then allow mutual interactions between the outcome and the covariates at the same time point.
For instance, we allow configurations of type $X_t=G(\eta_t,\eta_{t-1},\ldots)$ with $\left((\eta_t,\varepsilon_t)\right)_{t\in\Z}$ a sequence of i.i.d. random vectors such that $\eta_t$ is stochastically dependent of $\varepsilon_t$.
\subsubsection{Multinomial autoregressions}
In this part, we consider the finite state space $E=\{1,2,\ldots,N\}$ and some models satisfying
\begin{equation}\label{mainex2}
\P\left(Y_t=y\vert X_{t-1},Y_{t-1},\ldots\right)=H_y\left(X_{t-1},Y_{t-1},\ldots,Y_{t-p}\right),\quad y\in E,
\end{equation}
where $H_y$ is measurable mapping taking values in $(0,1)$ and such that $\sum_{y\in E}H_y=1$.
Such a framework include the multinomial logistic model for which for $y=1,\ldots,N-1$,
$$\log\left(H_y\left(X_{t-1},Y_{t-1},\ldots,Y_{t-p}\right)/H_N\left(X_{t-1},Y_{t-1},\ldots,Y_{t-p}\right)\right)=g\left(X_{t-1},Y_{t-2},\ldots,Y_{t-p}\right),$$
where $g$ is a measurable mapping that be linear but also quadratic or to exhibit more complex interactions between lag-values of the response and the covariates. See \citet{Truquet} for a discussion.
Here, we construct a random map by taking a uniformly distributed random variable $\varepsilon_t$ and setting
$$f\left(Y_{t-1},\ldots,Y_{t-p},X_{t-1},\varepsilon_t\right)=i \Leftrightarrow \sum_{1\leq y\leq i-1}H_{y,t}<\varepsilon_t\leq \sum_{1\leq y\leq i}H_{y,t},$$
where we set $H_{y,t}=H_y\left(X_{t-1},Y_{t-1},\ldots,Y_{t-p}\right)$ and use the convention $\sum_{y=1}^0 H_{y,t}=0$. 
We then obtain the following result.

\begin{prop}\label{seq1}
Suppose that Assumptions {\bf B2-B3} hold true. Then condition {\bf B1} is satisfied.
\end{prop}

\paragraph{Proof of Proposition \ref{seq1}}
Since $H_y$ is positive, we have 
$$\left\{\# F_t^{t+p}(E)=1\right\}\supset \cap_{i=0}^{p-1}\left\{\varepsilon_{t+i}\leq \inf_{y\in E^p} H_{1,t+i}\right\}.$$
Indeed, the intersection of these events leads to a value $1$ at time $t,t+1,\ldots,t+p-1$ and then to the same value at time $t+p$.
Now set $\mathcal{F}_t=\sigma\left((\varepsilon_s,X_s): s\leq t\right)$. Since
$$\P\left(\varepsilon_{t+i}\leq \inf_{y\in E^p} H_{1,t+i}\vert \mathcal{F}_{t+i-1}\right)=\inf_{y\in E^p}H_{1,t+i}>0,$$
one can apply Lemma \ref{ult2} to conclude that the intersection of the $p$ events is of positive probability. This concludes the proof.$\square$

\subsubsection{Ordinal time series}
Ordinal time series are typically constructed from a discretization of a continuous regression model. Here, we follow the presentation of \citet{Fokianos2003}, see in particular paragraph $3.2$ of that paper.
The state space of the process is $E=\{1,2,\ldots,N\}$. There is then a natural ordering on the state space and the autoregressive time series $(Y_t)_{t\in\Z}$ is defined by 
$$Y_t=i \Leftrightarrow c_{i-1}<g\left(Y_{t-1},\ldots,Y_{t-p},X_{t-1}\right)+\varepsilon_t\leq c_i,$$
where $-\infty=c_0<c_1<\cdots<c_{N-1}<c_N=\infty$, $g:E^p\times \R^d\rightarrow \R$ is a measurable function. 
The model writes in the form (\ref{mainex}) for a function $f$ that can be written as a linear combination of indicator sets.

\begin{prop}\label{ordinal}
If Assumptions {\bf B2-B3} are satisfied and the probability distribution of $\varepsilon_t$ has a density $f_{\varepsilon}$ which is positive on a given ray $(-\infty,r)$ or $(r,\infty)$,
then there exists a unique stationary and ergodic solution for (\ref{mainex}) and {\bf B1} is satisfied.
\end{prop}

\paragraph{Proof of Proposition \ref{ordinal}}

It is only necessary to check {\bf B1}. We show that such a condition is satisfied with $m=p$. 
Without loss of generality, assume that $f_{\varepsilon}$ is positive on the ray $(r,\infty)$. A similar argument applies when the ray is $(-\infty,r)$.
Indeed, it is clear that 
$$\left\{\# F_t^{t+p}(E)=1\right\}\supset \cap_{i=0}^{p-1}\left\{\varepsilon_{t+i}>c_{N-1}-\inf_{y\in E}g(y,X_{t+i-1})\right\}.$$
Indeed, if the event on the left occurs, we know that $p$ successive values $1$ will appear in the dynamic. As a consequence, at time $t$, the value of the iterated random maps will no more depend on the initial conditions. Since we have 
$$\P\left(\varepsilon_{t+i}>c_{N-1}-\inf_{y\in E}g(y,X_{t+i-1})\vert \mathcal{F}_{t+i-1}\right)=S\left(c_{N-1}-\inf_{y\in E}g(y,X_{t+i-1})\right),$$
where $S$ denotes the survival function of $\varepsilon_0$ which is positive, an application of Lemma \ref{ult2} shows that the intersection of the $p$ events has a positive probability and {\bf B1} is automatically satisfied.$\square$

\subsection{Dynamic multiple choice models}
Here we still consider the set $E=\{1,\ldots,N\}$ and we assume that  
$$Y_t=i\Leftrightarrow g_i\left(X_{t-1},Y_{t-1},\ldots,Y_{t-p}\right)+\varepsilon_{i,t}> g_j\left(X_{t-1},Y_{t-1},\ldots,Y_{t-p}\right)+\varepsilon_{j,t},\quad 1\leq j\neq i\leq N.$$
Such a model is called dynamic multiple choice model in theoretical economics and represents the successive choices of an agent that makes a decision after 
observing a set of covariates $X_{t-1}$.  See for instance \citet{linton} for semi-parametric inference in such models.

\begin{prop}\label{multiple}
Suppose that Assumptions {\bf B2-B3} hold true and that $\varepsilon_t$ has a distribution with a full support $\R^N$. Then {\bf B1} holds true.
\end{prop}

\paragraph{Proof of Proposition \ref{multiple}}

The proof is very similar than that of Proposition \ref{ordinal}. In particular, setting $g_{i,t}=g_i\left(X_{t-1},Y_{t-1},\ldots,Y_{t-p}\right)$, we have
$$\left\{\# F_t^{t+p}(E)=1\right\}\supset \cap_{j=0}^{p-1}\left\{\varepsilon_{i,t+j}>-g_{i,t+j}+\max_{\ell\neq i} \{g_{\ell,t+j}+\varepsilon_{\ell,t+\ell}\}\right\}$$
and the results follows from {\bf B2-B3}, the assumption of full support for the noise and Lemma \ref{ult2}.$\square$

\section{Sequential exogeneity and infinite dependence in autoregressive time series models}\label{5}

In this section, we first consider a general setup for infinite memory autoregressive processes with dependent inputs, in the spirit of Section \ref{4}.
Our results extend the models considered in \citet{Dvendange} when the inputs are independent. See also \citet{doukhantruquet} for similar results in the case of random fields.
We then use this result for deriving mixing properties of some autoregressive categorical time series or autoregressive time series of counts
which generalize some models found in the literature. Here the case of predetermined exogenous regressors is our main motivation.

\subsection{A general contraction argument}
We consider two Polish spaces $E$ and $G$ and we denote by $\Delta$ the metric used on $E$. Let $\left(\zeta_t\right)_{t\in\Z}$ be a stationary and ergodic process and $\overline{y}$ 
a reference point in $E$. Let $F:\mathcal{C}\times G\rightarrow E$ be a measurable mapping. We assume that for any nonnegative integer $k$,
the set $\mathcal{C}$ contains the subset of sequences $y$ of $E^{\N}$ such that $y_i=\overline{y}$ for $i\geq k$. Here $\overline{y}$ denotes a reference point in $E$.
We will also denote by $\overline{y}^{-}$ the element $(\overline{y},\overline{y},\ldots)$ of $E^{\N}$. 
In what follows, we set $\mathcal{F}_t=\sigma\left(\zeta_s:s\leq t\right)$. The two following assumptions will be used.

\begin{description}
\item[I1]
$\E \Delta\left(\overline{y},F(\overline{y}^{-},\zeta_0)\right)<\infty$.
\item[I2]
There exists a sequence $(a_i)_{i>0}$ of non-negative real numbers such that $a:=\sum_{i=1}^{\infty}a_i<1$ and for every $(y,y')\in \mathcal{C}^2$,
$$\E\left[\Delta\left(F(y,\zeta_1),F(y',\zeta_1)\right)\vert \mathcal{F}_0\right]\leq \sum_{i=1}^{\infty} a_i \Delta\left(y_i,y'_i\right)\mbox{ a.s.}.$$
\end{description}
The conditional contraction method generalizes the approach of \citet{DT} for finite-order autoregressive processes with dependent inputs.
Condition {\bf I2} has to be understood in term of regular conditional probability measure, i.e. there exists a regular version of the distribution of $\zeta_1$ conditionally on $\mathcal{F}_0$.

\begin{theo}\label{infinie}
Suppose that Assumptions {\bf I1-I2} hold true. There then exists a unique stationary and ergodic process $(Y_t)_{t\in\Z}$, non-anticipative, such that $\E \Delta(\overline{y},Y_0)<\infty$ and solution of
$$Y_t=F\left(Y_{t-1}^{-},\zeta_t\right),\quad t\in \Z,\mbox{ a.s}.$$
This unique solution has a Bernoulli shift representation, i.e.
$$Y_t=H\left(\zeta_t,\zeta_{t-1},\ldots\right)\mbox{ a.s.}$$
for some measurable function $H:\G^{\N}\rightarrow E$.
\end{theo}

\paragraph{Note.} By non-anticipative, we mean that $Y_t$ is $\mathcal{F}_t-$measurable for all $t\in \Z$. 
Let us mention that for any stationary and non-anticipative process $(Y_t)_{t\in\Z}$ taking values in $E$, the random variable $F\left(Y_{t-1}^{-},\zeta_t\right)$ is simply defined as the limit of $Z_{t,p}=F\left(Y_{t-1},\ldots,Y_{t-p},\overline{y}^{-},\zeta_1\right)$
when $p\rightarrow \infty$, in
$$\mathbb{L}^1_{\Delta}=\left\{Z:\Omega\rightarrow E\mbox{ measurable }: \E \Delta(\overline{y},Z)<\infty\right\}.$$
Indeed, since the process is non-anticipative, one can use {\bf I2} to get 
\begin{equation}\label{plusclair}
\E \Delta\left(Z_{t,p},Z_{t,q}\right)=\E\left[\E\left[\Delta(Z_{t,p},Z_{t,q})\right]\vert \mathcal{F}_{t-1}\right]\leq \sum_{i=p+1}^q a_i \E\left(\Delta(Y_1,\overline{y})\right),\quad p<q
\end{equation}
and $(Z_{t,p})_p$ is a Cauchy sequence which has then a limit.

\paragraph{Proof of Theorem \ref{infinie}}
Let $N$ be a positive integer. Define the process $\left(Y_{N,t}\right)_{t\in \Z}$ as follows. We set $Y_{N,t}=\overline{y}$ if $t\leq -N$ and for $t>N$, $Y_{N,t}$ is defined recursively by
$$Y_{N,t}=F\left(Y_{N,t-1}^{-},\zeta_t\right),\quad t\in\Z.$$
Since $Y_{N,t-1}^{-}$ is $\mathcal{F}_{t-1}-$measurable, we will use {\bf I2} in the following way.
Setting $M=\E d\left(\overline{y},F(\overline{y}^{-},\zeta_1)\right)$, we have for $t>-N$,
\begin{eqnarray*}
\E \Delta\left(Y_{N,t},\overline{y}\right)&\leq& M+\E \Delta\left(Y_{N,t},F(\overline{y}^{-},\zeta_t)\right)\\
&\leq & M+\sum_{i=1}^{\infty} a_i \E \Delta\left(Y_{N,t-i},\overline{y}\right).
\end{eqnarray*}
Using an induction argument, one can show that 
\begin{equation}\label{bound1}
\E \Delta\left(Y_{N,t},\overline{y}\right)\leq \frac{M}{1-a}.
\end{equation}
Set $C:=2M/(1-a)$. We are going to show that $\left(Y_{N,t}\right)_{N>-t}$ defines a Cauchy sequence in the space $\mathbb{L}_{\Delta}^{1}$. 
Let $N'>N$.
We are going to show that 
\begin{equation}\label{bound2}
\E \Delta\left(Y_{N,t},Y_{N',t}\right)\leq C\inf_{p\geq 1}\left\{a^{\frac{\max(0,t+N)}{p}}+\frac{S_{p+1}}{1-a}\right\},
\end{equation}
where $S_t=\sum_{j\geq t}a_j$ for $t\in\N^{*}$. Indeed, if we fix $N,p\geq 1$, we have $\Delta\left(Y_{N,t},Y_{N',t}\right)=0$ if $t\leq -N'$, 
$$\E \Delta\left(Y_{N,t},Y_{N',t}\right)=\E \Delta\left(\overline{y},Y_{N',t}\right)\leq C,\quad -N'<t\leq -N$$
and for $t>-N$, we have using {\bf I2} and the analogue of (\ref{plusclair}),
$$\E \Delta\left(Y_{N,t},Y_{N',t}\right)\leq \sum_{i=1}^p a_i \E \Delta\left(Y_{N,t-i},Y_{N',t-i}\right)+C S_{p+1}.$$
An induction argument for $t=-N+1,-N+2,\ldots$ leads to
$$\E \Delta\left(Y_{N,t},Y_{N',t}\right)\leq C\left\{a^{\frac{\max(0,t+N)}{p}}+\frac{S_{p+1}}{1-a}\right\}$$
and then to (\ref{bound2}).  
For any $t\in \Z$, let us denote by $Y_t$ the limit of $Y_{N,t}$. Now let $\mu$ be the probability distribution on $\left(\zeta_{-i}\right)_{i\geq 1}$. 
Since there exists a measurable function $H_N: G^{\N}\rightarrow E$ such that $Y_{N,t}=H_N\left(\zeta_t,\zeta_{t-1},\ldots\right)$ and the sequence $(H_N)_{N\in \N}$ is a Cauchy sequence 
in the space $\mathcal{G}$ of measurable functions $H:G^{\N}\rightarrow E$ satisfying $\int \Delta(H(x),\overline{y})d\mu(x)>\infty$ and endowed with the distance $(H,H')\mapsto  \int \Delta\left(H(x),H'(x)\right)d\mu(x)$,
one can easily find an element $H$ of $\mathcal{G}$ such that $Y_t=H\left(\zeta_t,\zeta_{t-1},\ldots\right)$ a.s.
This representation also entails stationary and ergodicity of the process $(Y_t)_{t\in\Z}$.
Moreover, for $\epsilon>0$, if $p$ is large enough, we have for any integer $N\geq 1$,
$$\E \Delta\left(Y_{N,t},F\left(Y_{N,t-1},\ldots,Y_{N,t-p},\overline{y}^{-},\zeta_t\right)\right)+\E \Delta\left(Y_t,F\left(Y_{t-1},\ldots,Y_{t-p},\overline{y}^{-},\zeta_t\right)\right)<\epsilon.$$
Moreover, in $\mathbb{L}_{\Delta}^1$,
$$\lim_{N\rightarrow \infty}F\left(Y_{N,t-1},\ldots,Y_{N,t-p},\overline{y}^{-},\zeta_t\right)=F\left(Y_{t-1},\ldots,Y_{t-p},\overline{y}^{-},\zeta_t\right).$$
This shows that 
$$\lim_{N\rightarrow \infty}F\left(Y_{N,t-1}^{-},\zeta_t\right)=F\left(Y_{t-1}^{-},\zeta_t\right)$$
in $\mathbb{L}_{\Delta}^1$. Then $(Y_t)_{t\in\Z}$ is indeed solution of the infinite memory autoregression equations.
Uniqueness follows from standard arguments using {\bf I2}. details are omitted.$\square$

We now particularize our results to the case $G=(\R^d)^{\N}\times \R^k$ and $\zeta_t=\left(X_{t-1}^{-},\varepsilon_t\right)$.
The regression function $F$ is now defined on $\mathcal{C}\times \mathcal{D}\times\R^k$ where $\P\left(X_0^{-}\in\mathcal{D}\right)=1$ and $\mathcal{D}$ contains all sequences in $(\R^d)^{\N}$ 
that vanishes for large indices. We are interested here in stationary solutions of
\begin{equation}\label{jenaim}
Y_t=F\left(Y_{t-1}^{-},X_{t-1}^{-},\varepsilon_t\right),\quad t\in\Z
\end{equation}
and we use the following specific condition. In what follows, $\vert\cdot\vert$ denotes an arbitrary norm on $\R^d$.
\begin{description}
\item[I2'] 
There exist two sequences of non-negative real numbers $(a_i)_{i\geq 1}$ and $(b_j)_{j\geq 1}$ such that $\sum_{i\geq 1} a_i<1$, $\sum_{j\geq 1}b_j<\infty$ and if $(y,y',x,x')\in\mathcal{C}^2\times \mathcal{D}^2$,
$$\E \Delta\left(F(y,x,\zeta_1),F(y',x',\zeta_1)\right)\leq \sum_{i\geq 1}a_i \Delta(y_i,y'_i)+\sum_{j\geq 1} b_j \vert x_j-x'_j\vert.$$
\end{description}
Here, $\mathcal{F}_t=\sigma\left((X_s,\varepsilon_s): s\leq t\right)$. 

\begin{theo}\label{final}
Suppose that Assumptions {\bf B2-B3} and {I1-I2'} hold true. Suppose furthermore that $\E \vert X_0\vert<\infty$.
\begin{enumerate}
\item
There then exists a unique stationary and ergodic 
$(Y_t)_{t\in\Z}$, non-anticipative, such that $\E d(\overline{y},Y_0)<\infty$ and solution of
$$Y_t=F\left(Y_{t-1}^{-},X_{t-1}^{-},\varepsilon_t\right),\quad t\in \Z,\mbox{ a.s}.$$
This unique solution has a Bernoulli shift representation, i.e.
$$Y_t=H\left(X_{t-1},\varepsilon_t,X_{t-2},\varepsilon_{t-1},\ldots\right)\mbox{ a.s.}$$
\item
For any $0<r<n$, there exists a process $(Y'_t)_{t\in\Z}$ such that for $t\geq r+1$, $Y'_t$ is measurable with respect to $\sigma\left((X_j,\varepsilon_{j+1}): j\geq r\right)$ and 
$$\E \Delta(Y_t,Y'_t)\leq L\inf_{p\geq 1}\left\{a^{\frac{t-r}{p}}+S_{p+1}+\sum_{j=0}^{t-r+1}a^{j/p}T_{t-r+1-j}\right\},$$
where $L$ is a positive constant, $S_{p+1}=\sum_{i\geq p+1}a_i$ and $T_s=\sum_{j\geq s}b_j$ for $s\in\N$.
\end{enumerate}
\end{theo}

\paragraph{Proof of Theorem \ref{final}}
The first part about existence of a stationary solution is a consequence of Theorem \ref{infinie}. Note that {\bf B3-I2'} entails {\bf I2}.

Now, let $0<r<n$. Let us define $(Y'_t)_{t\in\Z}$ the process defined by $Y'_t=\overline{y}$ if $t\leq r$ and 
$${Y'}_t=F\left({Y'}_{t-1}^{-},X_{t-1},\ldots,X_r,0,\ldots,\varepsilon_t\right),\quad t\geq r+1.$$
We set $M=\E \Delta\left(F(\overline{y}^{-},X_0^{-},\varepsilon_1),\overline{y}\right)$ and $D_2=\E\vert X_0\vert$.
Note that is $t>r$, using the triangular inequality
\begin{eqnarray*}
h_t:=\Delta\left(Y'_t,\overline{y}\right)&\leq & \Delta\left(Y_t',F(\overline{y}^{-},X_{t-1},\ldots,X_r,0,\ldots,\varepsilon_t)\right)\\
&+& \Delta\left(F(\overline{y}^{-},X_{t-1},\ldots,X_r,0,\ldots,\varepsilon_t),F(\overline{y}^{-},X_{t-1}^{-},\varepsilon_t)\right)\\
&+& \Delta\left(F(\overline{y}^{-},X_{t-1}^{-},\varepsilon_t),\overline{y}\right)
\end{eqnarray*}
our assumptions guarantee that 
$$h_t\leq \sum_{i\geq 1}a_i h_{t-i}+\sum_{j\geq t-r+1}b_j \E\vert X_0\vert+M.$$
Since $h_t=0$ for $t\leq r$, there exists a constant $C_1>0$ such that $\sup_{t\in\Z}h_t\leq C_1$. Set $D_1=C_1+\E \Delta(Y_0,\overline{y})$.
Now if $t=r+s$, we have for $s\geq 1$,
\begin{eqnarray*}
u_s:=\E \Delta\left(Y_t,Y_t'\right)&\leq& \E \Delta\left(Y_t,F({Y'}_{t-1}^{-},X_{t-1}^{-},\varepsilon_t)\right)\\&+&\E \Delta\left(F({Y'}_{t-1}^{-},X_{t-1}^{-},\varepsilon_t),F({Y'}_{t-1}^{-},X_{t-1},\ldots,X_r,0,\ldots,\varepsilon_t)\right)\\
&\leq& \sum_{i\geq 1}a_i u_{s-i}+\sum_{j\geq s+1}b_j \E\vert X_0\vert\\
&\leq& \sum_{i=1}^pa_i u_{s-i}+D_1 S_{p+1}+D_2 T_{s+1}.
\end{eqnarray*}
Note that for $s\leq 0$, $u_s\leq C:=\E \Delta(Y_0,\overline{y})$. The result easily follows from an application of Lemma \ref{ult3}.$\square$

\subsection{Application to some binary time series} 
Here, we consider an example of binary process for which the covariates only predetermine the response. We set $E=\{0,1\}$ and $\Delta$ is the discrete metric, i.e. 
$\Delta(y,y')=\mathds{1}_{y\neq y'}$ for $y,y'\in E$.
Conditionally on $\left(Y_{t-1}^{-},X_{t-1}^{-}\right)$, we assume that $Y_t$ follows a Bernoulli distribution with random parameter $F\left(\lambda_t\right)$ and
\begin{equation}\label{dyndyn}
\lambda_t=f\left(Y_{t-1}^{-},X_{t-1}^{-}\right).
\end{equation}
Here $F$ is a c.d.f. such that the logistic or the standard Gaussian c.d.f.
Once again, this setup corresponds to some models introduced in \citet{rydberg} or \citet{Fok1}.
A natural way for considering stochastic iterations is to use the inverse of a c.d.f. associated to the Bernoulli distribution:
\begin{equation}\label{dyndyn+}
Y_t=\mathds{1}_{\{\varepsilon_t>1-p_t\}},\quad \varepsilon_t\sim\mathcal{U}([0,1]).
\end{equation} 

We remind that $\zeta_t=(X_{t-1},\varepsilon_t)$ and $V_t=(\zeta_t,Y_t)$.

\begin{cor}\label{binarypred}
Suppose that Assumptions {\bf B2-B3} hold true. Suppose furthermore that $F$ is Lipschitz with Lipschitz constant $L_F$ and that the exist two sequences 
$(a_i')_{i\geq 1}$ and $(b'_j)_{j\geq 1}$ of non-negative real numbers
\begin{enumerate}
\item
$L_F\sum_{i\geq 1}a_i'<1$ and $\sum_{j\geq 1}b'_j<\infty$.
\item
If $(y,y',x,,x')\in\mathcal{C}^2\times \mathcal{D}^2$,
$$\vert f(y,x)-f(y',x')\vert\leq \sum_{i\geq 1}a'_i\mathds{1}_{y_i\neq y_i'}+\sum_{j\geq 1} b'_j \vert x_j-x'_j\vert.$$
\end{enumerate}
There then exists a unique non-anticipative stationary and ergodic solution for (\ref{dyndyn})-(\ref{dyndyn+}). 
Moreover, for any $1<r<n$, we have the bound
$$\alpha_V(n)\leq \alpha_{\zeta}(r+1)+L\sum_{t\geq n}\inf_{p\geq 1}\left\{a^{\frac{t-r}{p}}+S_{p+1}+\sum_{j=0}^{t-r+1}a^{j/p}T_{t-r+1-j}\right\},$$
with $a_i=a_i'L_f$, $b_j=b_j'L_f$, $S_{p+1}=\sum_{i\geq p+1}a_i$ and $T_s=\sum_{j\geq s}b_j$ for $s\in \N$.
\end{cor} 

\paragraph{Note.}
When $\lambda_t=\beta \lambda_{t-1}+\kappa Y_{t-1}+\delta' X_{t-1}$, we recover the model of \citet{Fok1}. If $F(s)=(1+e^{-s})^{-1}$ (logistic case), we have $L_F=1/4$ and we get representation (\ref{dyndyn}) with $f(y,x)=\sum_{i\geq 1}\beta^{i-1}\left(\kappa y_i+\delta'x_i\right)$. Our result apply as soon as $\vert \beta\vert+\vert \gamma\vert/4<1$, which is the same stationary condition than that of \citet{Fok1}.
The coefficients $a_i$ and $b_j$ have exponential decays and we $\alpha_V(n)=O\left(\alpha_{\zeta}([n/2])\right)$ as soon as $\alpha_{\zeta}$ has a sub exponential decay. Note that the stationary condition is much more restrictive than in the case of strict exogeneity. In the latter case, $\vert\beta\vert<1$ is the needed condition, as explained in the Note just of Theorem \ref{mixing2}.

\paragraph{Proof of Corollary \ref{binarypred}}

The result is a simple consequence of Theorem \ref{final} with the discrete metric $\Delta(y,y')=\mathds{1}_{y\neq y'}$. In particular, one can use the bound 
$$\P\left(\mathds{1}_{\{\varepsilon_1>1-q\}}\neq \mathds{1}_{\{\varepsilon_1>1-q\}}\right)=\vert q-q'\vert, \quad q,q'\in [0,1],$$
to check {\bf I2'}.$\square$

\subsection{Application to INGARCH models}\label{6}
For this last class of examples, we consider $E=\R$ and $\Delta$ the distance defined by the absolute value.
The model is defined by 
\begin{equation}\label{integer}
Y_t=g\circ F^{-1}_{\lambda_t}(\varepsilon_t),\quad \lambda_t=f\left(Y_{t-1}^{-},X_{t-1}^{-}\right),
\end{equation}
where $F_{\lambda}$ denotes the c.d.f. of a probability distribution on $\N$, that depends on a real-valued parameter $\lambda$, $f,g$ are measurable functions
such that $g$ is invertible and $\varepsilon_t$ is uniformly distributed over the interval $[0,1]$.
A standard example concerns INGARCH processes, for which $F_{\lambda}$ is the c.d.f. of the Poisson distribution of parameter $\lambda$ and $g$ is the identity function. See for instance \citet{Doukhan2012c} and the references therein.
\citet{davis} considered many examples of integer-valued autoregressive models using the representation \ref{integer}, with $g(y)=y$ and $\lambda$ being the mean of $F_{\lambda}^{-1}(U_1)$.
Another example concerns the log-INGARCH model of \citet{fokianos2011log} for which $g(y)=\log(1+y)$ and $F_{\lambda}$ is the Poisson distribution with parameter $\exp(\lambda)$. In this case the function $f$ can take values of arbitrary sign. 
Let us mention that these models are generally presented under the observation-driven form, i.e.
\begin{equation}\label{usual}
\lambda_t=\beta \lambda_{t-1}+\kappa Y_{t-1}+\delta' X_{t-1}.
\end{equation}
and mainly without exogenous regressors (i.e. $\delta=0$). However, some recent references such as \citet{PARX} or \citet{DT3} studied existence of stationary solutions with the problem of covariates inclusion.  Note that under the condition $\vert \beta\vert<1$, which is always assumed, and an integrability condition on $X_0$, it is not difficult to show that stationary solutions of (\ref{usual}) coincide to stationary solutions of (\ref{integer}) with a function $f$ written as a series.

Our main result is the following. For the mixing coefficients, we distinguish the two cases $g(y)=y$ and $g(y)=\log(1+y)$ in the case of the log-INGARCH model. For simplicity of notations, for the coefficients $(a_i)_{i\geq 1}$ and $(b_j)_{j\geq 1}$ introduced below, we set
$$\omega_{n,t,r}=\inf_{p\geq 1}\left\{a^{\frac{t-r}{p}}+S_{p+1}+\sum_{j=0}^{t-r+1}a^{j/p}T_{t-r+1-j}\right\},$$
with $S_{p+1}=\sum_{i\geq p+1}a_i$ and $T_s=\sum_{j\geq s}b_j$ for $s\in \N^*$.
\begin{cor}\label{cestfini}
Suppose that Assumptions {\bf B2-B3} hold true and $\E\vert X_0\vert<\infty$. Suppose furthermore that 
\begin{enumerate}
\item
For any value of $\lambda$, $g\circ F_{\lambda}^{-1}(\varepsilon_1)$ and $f\left(\overline{y}^{-},X_0^{-1}\right)$ are integrable and 
$$\E\left\vert g\circ F_{\lambda}^{-1}(\varepsilon_1)- g\circ F_{\lambda'}^{-1}(\varepsilon_1)\right\vert \leq \vert \lambda-\lambda'\vert.$$
\item
There exist two summable sequences of non-negative real numbers $(a_i)_{i\geq 1}$ and $(b_j)_{j\geq 1}$ such that $\sum_{i\geq 1} a_1<1$ and 
for $(y,y',x, x')\in\mathcal{C}^2\times \mathcal{D}^2$,
$$\vert f(y,x)-f(y',x')\vert\leq \sum_{i\geq 1}a_i\vert y_i-y_i'\vert+\sum_{j\geq 1} b_j \vert x_j-x'_j\vert.$$
\end{enumerate}
There then exists a unique non-anticipative stationary and ergodic solution for (\ref{integer}).
 
Moreover, for any $1<r<n$, we have the following bounds.
\begin{enumerate}
\item
If $g(y)=y$,
$$\alpha_V(n)\leq \alpha_{\zeta}(r+1)+L\sum_{t\geq n}\omega_{n,t,r}.$$
\item
If $g(y)=\log(1+y)$, $F_{\lambda}$ corresponds to the Poisson distribution with parameter $\exp(\lambda)$ and $\E\exp\left(K \vert X_0\vert\right)<\infty$ for any $K\geq 1$, then for any $K\geq 1$, there exists a constant $L>0$ such that  
$$\alpha_V(n)\leq \alpha_{\zeta}(r+1)+L\sum_{t\geq n}\omega_{n,t,r}^{K/(K+1)}.$$ 
\end{enumerate}
\end{cor}

\paragraph{Note.}
For INGARCH models as in (\ref{usual}), one can apply the previous result as soon as $\vert \beta\vert +\vert\kappa\vert<1$. 
In this case and as for the previous example of binary process, we obtain $\alpha_V(n)=O\left(\alpha_{\zeta}([n/2])\right)$ as soon as the mixing coefficients of $\zeta$ have a sub exponential decay.
Note also that the first assumption is satisfied for both the Poisson INGARCH and the log-INGARCH model. See for instance \citet{DT2}, Section $2.3$.
Let us mention that for the log-INGARCH model, we also obtain a control of the strong mixing coefficients of $\left(g^{-1}(Y_t)\right)_{t\in \Z}$ which is the integer-valued process. Indeed, the mixing coefficients are invariant under a bijective transformation of the process.

\paragraph{Proof of Corollary \ref{cestfini}}

We proceed as for the proof of Corollary \ref{binarypred}, using Theorem \ref{final}. To this end, we use the metric $\Delta(y,y')=\vert y-y\vert$ on $E=\R$.
Existence of a unique stationary solution is straightforward to get. Let us discuss, the bound on the mixing coefficients.
\begin{enumerate}
\item
If $g(y)=y$, the single missing argument is a bound on $\P\left(Y_t\neq Y'_t\right)$ obtained form the the bound $\E\Delta(Y_t,Y'_t)$ given in Theorem \ref{final}. It is simply necessary to note that $\P\left(Y_t\neq Y'_t\right)\leq \E\Delta(Y_t,Y_t')$ because $Y_t$ and $Y'_t$ are integer-valued.
\item
If now $g(y)=\log(1+y)$, only $Z_t=g^{-1}(Y_t)$ and $Z_t'=g^{-1}(Y_t')$ are integer valued. We have for any real number $C>0$ and $K>0$,
\begin{eqnarray*}
\P\left(Y_t\neq Y'_t\right)&=&\P\left(Z_t\neq Z'_t\right)\\
&\leq& \P\left(Z_t\neq Z'_t, Z_t\leq C,Z'_t\leq C\right)+\P\left(Z_t>C\right)+\P\left(Z'_t>C\right)\\
&\leq& \E\left[\Delta\left(Z_t,Z'_t\right)\mathds{1}_{\max(Z_t,Z'_t)\leq C}\right]+\frac{\E Z_0^K}{C^K}+\frac{\sup_{t\geq r}\E{{Z'}_t^K}}{C^K}\\
&\leq& (1+C)\E\left[\Delta\left(Y_t,Y'_t\right)\right]+\frac{\E Z_0^K}{C^K}+\frac{\sup_{t\geq r}\E{{Z'}_t^K}}{C^K}.
\end{eqnarray*}
The last inequality is obtained from an application of the mean value theorem, i.e. $\vert \log(1+z)-\log(1+z')\vert\geq (1+C)^{-1}\vert z-z'\vert$ when 
$\max(z,z')\leq C$.
One can then choose $C\sim \E\left[\Delta\left(Y_t,Y'_t\right)\right]^{-K/(K+1)}$ to finally bound $\P\left(Y_t\neq Y'_t\right)$ by $\E\left[\Delta\left(Y_t,Y'_t\right)\right]^{-K/(K+1)}$ and apply Theorem \ref{final}. The single point to justify is the moment condition
$$\E Z_0^K+ \sum_{t\geq r}\E {Z'}_t^K<\infty.$$
We will not detail all the arguments for this but we simply mention that from our assumptions and a convexity argument, we have for any $K\geq 1$,
$$\E\exp\left(K\sum_{j\geq 1}b_j \vert X_{t-j}\vert\right)<\infty\mbox{ and } \E\exp\left(K\left\vert f\left(\overline{y}^{-},X_{t-1}^{-}\right)\right\vert\right)<\infty.$$
Moreover, using convexity arguments similar to \citet{DT3}, proof of Proposition $4$, one can show $g^{-1}(Y_t)$ and $g^{-1}(Y'_t)$ have moments of any order.
Details are omitted.$\square$
\end{enumerate}
\section{Appendix}\label{6}

\begin{lem}\label{ult}
Suppose that $\left(\kappa_t\right)_{t\in \Z}$ is a stationary process taking values in the interval $[0,1]$. Set $\rho=\E\kappa_1<1$. Then 
\begin{equation}\label{upbd}
\E\left(\kappa_1\cdots \kappa_s\right)\leq \inf_{1\leq j\leq s-1}\left\{\rho^{[s/j]}+\frac{\alpha_{\kappa}(j)}{1-\rho}\right\}.
\end{equation}
\end{lem}

\paragraph{Proof of Lemma \ref{ult}}
Setting $u_s=\E\left(\kappa_1\cdots \kappa_s\right)$ and using that all the variables take their values in the unit interval, we have
\begin{eqnarray*}
u_s&\leq &\E\left(\kappa_1\cdots\kappa_{s-j}\kappa_s\right)\\
&\leq& \rho u_{s-j}+\left\vert \c\left(\kappa_s,\kappa_1\cdots\kappa_{s-j}\right)\right\vert\\
&\leq & \rho u_{s-j}+4\alpha_{\kappa}(j),
\end{eqnarray*}
where the last inequality follows from the covariance inequality given in \citet{Doukhan(1994)}, Lemma $3$. 
Since $\rho<1$ and $u_k\leq 1$ for $1\leq k\leq s$, one can iterate the previous bound to get 
$$u_s\leq \rho^{[s/j]}+\frac{\alpha_{\kappa}(j)}{1-\rho}.$$
Taking the infinimum over $j$ leads to the proposed bound.$\square$

\begin{lem}\label{ult2}
Let $n\in \N$. Suppose that $\left(\mathcal{F}_t\right)_{0\leq t\leq n}$ is a filtration and $A_t\in\mathcal{F}_t$,  $1\leq t\leq n$, are some events such that $\P\left(A_t\vert\mathcal{F}_{t-1}\right)>0$ a.s.
Then $\P\left(\cap_{t=1}^n A_t\right)>0$.
\end{lem}

\paragraph{Proof of Lemma \ref{ult2}}
We show the result by induction on $n$.
From our assumptions, we know that $\P(A_1)>0$. Then the result is true for $n=1$. Suppose that the result holds true for $n=j$.
We have 
$$\P\left(\cap_{\ell=1}^{j+1}A_{\ell}\right)=\P\left(\P\left(A_{j+1}\vert \mathcal{F}_j\right)\mathds{1}_{\cap_{\ell=1}^jA_{\ell}}\right).$$
If this probability equals to $0$, our assumption entails that $\P\left(\cap_{\ell=1}^j A_{\ell}\right)=0$ which contradicts the induction assumption.
The result then follows.$\square$

\begin{lem}\label{ult3}
Assume that $(u_t)_{t\in\Z}$ is a sequence of real numbers such that $u_t\leq C$ for $t\leq 0$ and $u_t\leq \sum_{i=1}^p a_i u_{t-i}+v_t$ for $t>0$ where $(v_t)_{t>0}$ is a non-increasing sequence and $(a_i)_{i\geq 1}$ is a sequence of non-negative real numbers such that $\sum_{i=1}^{\infty} a_i<1$. We then get 
$$u_t\leq C a^{\frac{t\vee 0}{p}}+\sum_{j\geq 0}a^{j/p}v_{(t-j)\vee 0}.$$
\end{lem}

\paragraph{Proof of Lemma \ref{ult3}}

The announced bound is easily checked when $t\leq 0$. If the bound is satisfied up to time $t\geq 0$, we note that 
\begin{eqnarray*}
u_{t+1}&\leq & C\sum_{i=1}^p a_i a^{\frac{(t+1-i)\vee 0}{p}}+\sum_{j\geq 0}a^{j/p}\sum_{i=1}^p a_i v_{(t+1-i-j)\vee 0}+v_{t+1}\\
&\leq& C a a^{\frac{(t+1-p)\vee 0}{p}}+\sum_{j\geq 0}a^{(j+p)/p}v_{(t+1-p-j)\vee 0}+v_{t+1}\\
&\leq & C a^{\frac{t+1}{p}}+\sum_{j\geq 0}a^{j/p}v_{(t+1-j)\vee 0}.
\end{eqnarray*}
The result then follows from an induction argument.$\square$
\bibliographystyle{plainnat}
\bibliography{biblio}

\begin{thebibliography}{32}
\providecommand{\natexlab}[1]{#1}
\providecommand{\url}[1]{\texttt{#1}}
\expandafter\ifx\csname urlstyle\endcsname\relax
  \providecommand{\doi}[1]{doi: #1}\else
  \providecommand{\doi}{doi: \begingroup \urlstyle{rm}\Url}\fi

\bibitem[Agosto et~al.(2016)Agosto, Cavaliere, Kristensen, and Rahbek]{PARX}
Arianna Agosto, Giuseppe Cavaliere, Dennis Kristensen, and Anders Rahbek.
\newblock Modeling corporate defaults: Poisson autoregressions with exogenous
  covariates (parx).
\newblock \emph{Journal of Empirical Finance}, 38:\penalty0 640--663, 2016.

\bibitem[Chamberlain(1982)]{Chamb}
Gary Chamberlain.
\newblock The general equivalence of granger and sims causality.
\newblock \emph{Econometrica: Journal of the Econometric Society}, pages
  569--581, 1982.

\bibitem[Chazottes et~al.(2020)Chazottes, Gallo, and Takahashi]{Chaz}
J-R Chazottes, S.~Gallo, and D.~Takahashi.
\newblock Optimal gaussian concentration bounds for stochastic chains of
  unbounded memory.
\newblock \emph{arXiv preprint arXiv:2001.06633}, 2020.

\bibitem[Davis and Liu(2016)]{davis}
Richard~A Davis and Heng Liu.
\newblock Theory and inference for a class of nonlinear models with application
  to time series of counts.
\newblock \emph{Statistica Sinica}, pages 1673--1707, 2016.

\bibitem[de~Jong and Woutersen(2011)]{deJ}
Robert~M de~Jong and Tiemen Woutersen.
\newblock Dynamic time series binary choice.
\newblock \emph{Econometric Theory}, 27\penalty0 (4):\penalty0 673--702, 2011.

\bibitem[Debaly and Truquet(2019)]{DT2}
Zinsou~Max Debaly and Lionel Truquet.
\newblock Stationarity and moment properties of some multivariate count
  autoregressions.
\newblock \emph{arXiv preprint arXiv:1909.11392}, 2019.

\bibitem[Debaly and Truquet(2021{\natexlab{a}})]{DT}
Zinsou~Max Debaly and Lionel Truquet.
\newblock Iteration of dependent random maps and exogeneity in nonlinear
  dynamics.
\newblock \emph{Econometric Theory}, pages 1--38, 2021{\natexlab{a}}.
\newblock \doi{10.1017/S0266466620000559}.

\bibitem[Debaly and Truquet(2021{\natexlab{b}})]{DT3}
Zinsou~Max Debaly and Lionel Truquet.
\newblock Multivariate time series models for mixed data.
\newblock \emph{arXiv preprint arXiv:2104.01067}, 2021{\natexlab{b}}.

\bibitem[Doukhan(1994)]{Doukhan(1994)}
P.~Doukhan.
\newblock \emph{Mixing: properties and examples}.
\newblock Number~85 in Lecture Notes in Statistics. Springer-Verlag, New York,
  1994.

\bibitem[Doukhan and Truquet(2007)]{doukhantruquet}
Paul Doukhan and Lionel Truquet.
\newblock Weakly dependent random fields with infinite interactions-paru sous
  le titre" a fixed point approach to model random fields".
\newblock \emph{ALEA: Latin American Journal of Probability and Mathematical
  Statistics}, 3:\penalty0 111--132, 2007.

\bibitem[Doukhan and Wintenberger(2008)]{Dvendange}
Paul Doukhan and Olivier Wintenberger.
\newblock Weakly dependent chains with infinite memory.
\newblock \emph{Stochastic Processes and their Applications}, 118\penalty0
  (11):\penalty0 1997--2013, 2008.

\bibitem[Doukhan et~al.(2012{\natexlab{a}})Doukhan, Fokianos, and
  Li]{doukhan2012weak}
Paul Doukhan, Konstantinos Fokianos, and Xiaoyin Li.
\newblock On weak dependence conditions: The case of discrete valued processes.
\newblock \emph{Statistics \& Probability Letters}, 82\penalty0 (11):\penalty0
  1941--1948, 2012{\natexlab{a}}.

\bibitem[Doukhan et~al.(2012{\natexlab{b}})Doukhan, Fokianos, and
  Tj{\o}stheim]{Doukhan2012c}
Paul Doukhan, Konstantinos Fokianos, and Dag Tj{\o}stheim.
\newblock On weak dependence conditions for {P}oisson autoregressions.
\newblock \emph{Statist. Probab. Lett.}, 82:\penalty0 942--948,
  2012{\natexlab{b}}.

\bibitem[Doukhan et~al.(2020)Doukhan, Neumann, and Truquet]{DNT}
Paul Doukhan, Michael~H Neumann, and Lionel Truquet.
\newblock Stationarity and ergodic properties for some observation-driven
  models in random environments.
\newblock \emph{arXiv preprint arXiv:2007.07623}, 2020.

\bibitem[Fern{\'a}ndez-Val and Weidner(2016)]{Fern}
Iv{\'a}n Fern{\'a}ndez-Val and Martin Weidner.
\newblock Individual and time effects in nonlinear panel models with large n,
  t.
\newblock \emph{Journal of Econometrics}, 192\penalty0 (1):\penalty0 291--312,
  2016.

\bibitem[Fokianos and Kedem(2003)]{Fokianos2003}
Konstantinos Fokianos and Benjamin Kedem.
\newblock Regression theory for categorical time series.
\newblock \emph{Statist. Sci.}, 18:\penalty0 357--376, 2003.
\newblock ISSN 0883-4237.
\newblock \doi{10.1214/ss/1076102425}.
\newblock URL \url{http://dx.doi.org/10.1214/ss/1076102425}.

\bibitem[Fokianos and Tj{\o}stheim(2011)]{fokianos2011log}
Konstantinos Fokianos and Dag Tj{\o}stheim.
\newblock Log-linear poisson autoregression.
\newblock \emph{Journal of Multivariate Analysis}, 102\penalty0 (3):\penalty0
  563--578, 2011.

\bibitem[Fokianos and Truquet(2018)]{Truquet}
Konstantinos Fokianos and Lionel Truquet.
\newblock On categorical time series models with covariates.
\newblock \emph{Stochastic Processes and their Applications}, 2018.

\bibitem[Herv{\'e} and Ledoux(2021)]{HL}
Lo{\"\i}c Herv{\'e} and James Ledoux.
\newblock Asymptotic of products of markov kernels. application to
  deterministic and random forward/backward products.
\newblock \emph{Statistics \& Probability Letters}, 179:\penalty0 109204, 2021.

\bibitem[Hsiao(2014)]{Hsiao}
Cheng Hsiao.
\newblock \emph{Analysis of panel data}.
\newblock Number~54. Cambridge university press, 2014.

\bibitem[Kauppi and Saikkonen(2008)]{kauppi}
Heikki Kauppi and Pentti Saikkonen.
\newblock Predicting us recessions with dynamic binary response models.
\newblock \emph{The Review of Economics and Statistics}, 90\penalty0
  (4):\penalty0 777--791, 2008.

\bibitem[Kifer(1996)]{Kifer}
Yuri Kifer.
\newblock Perron-frobenius theorem, large deviations, and random perturbations
  in random environments.
\newblock \emph{Mathematische Zeitschrift}, 222\penalty0 (4):\penalty0
  677--698, 1996.

\bibitem[Moysiadis and Fokianos(2014)]{Fok1}
T.~Moysiadis and K.~Fokianos.
\newblock On binary and categorical time series models with feedback.
\newblock \emph{J. Multivariate Anal.}, 131:\penalty0 209--228, 2014.

\bibitem[Neumann(2011)]{Neumann2011}
M.~Neumann.
\newblock Absolute regularity and ergodicity of poisson count processes.
\newblock \emph{Bernoulli}, 17:\penalty0 1268--1284, 2011.

\bibitem[Park et~al.(2017)Park, Simar, and Zelenyuk]{park}
Byeong~U Park, L{\'e}opold Simar, and Valentin Zelenyuk.
\newblock Nonparametric estimation of dynamic discrete choice models for time
  series data.
\newblock \emph{Computational Statistics \& Data Analysis}, 108:\penalty0
  97--120, 2017.

\bibitem[Park and Phillips(2000)]{Ph}
Joon~Y Park and Peter~CB Phillips.
\newblock Nonstationary binary choice.
\newblock \emph{Econometrica}, 68\penalty0 (5):\penalty0 1249--1280, 2000.

\bibitem[Rocha and Cribari-Neto(2009)]{beta}
Andr{\'e}a~V Rocha and Francisco Cribari-Neto.
\newblock Beta autoregressive moving average models.
\newblock \emph{Test}, 18\penalty0 (3):\penalty0 529--545, 2009.

\bibitem[Russell and Engle(2005)]{russell}
Jeffrey~R Russell and Robert~F Engle.
\newblock A discrete-state continuous-time model of financial transactions
  prices and times: The autoregressive conditional multinomial--autoregressive
  conditional duration model.
\newblock \emph{Journal of Business \& Economic Statistics}, 23\penalty0
  (2):\penalty0 166--180, 2005.

\bibitem[Rydberg and Shephard(2003)]{rydberg}
Tina~Hviid Rydberg and Neil Shephard.
\newblock Dynamics of trade-by-trade price movements: decomposition and models.
\newblock \emph{Journal of Financial Econometrics}, 1\penalty0 (1):\penalty0
  2--25, 2003.

\bibitem[Srisuma and Linton(2012)]{linton}
Sorawoot Srisuma and Oliver Linton.
\newblock Semiparametric estimation of markov decision processes with
  continuous state space.
\newblock \emph{Journal of Econometrics}, 166\penalty0 (2):\penalty0 320--341,
  2012.

\bibitem[Truquet(2020)]{Truquet2}
Lionel Truquet.
\newblock Coupling and perturbation techniques for categorical time series.
\newblock \emph{Bernoulli}, 26\penalty0 (4):\penalty0 3249--3279, 2020.

\bibitem[Wei{\ss}(2018)]{Weiss}
Christian~H Wei{\ss}.
\newblock \emph{An introduction to discrete-valued time series}.
\newblock John Wiley \& Sons, 2018.

\end{thebibliography}

\end{document}